\documentclass[12pt]{article}

\usepackage{amsmath, amsthm}
\usepackage{amsfonts}
\usepackage{amssymb}
\usepackage[usenames,dvipsnames]{color}%\usepackage[dvips]{color}
\usepackage{graphicx}
\usepackage{hyperref}
\usepackage{bm}
\usepackage{bbm}
\usepackage{listings}
\usepackage{geometry}
\usepackage{lipsum} 

\usepackage{ulem}

\usepackage{color,calc}
\definecolor{junglegreen}{rgb}{0.16, 0.67, 0.53}

\makeatletter
\definecolor{shade}{gray}{0.8}
        {
          \raggedright
        \setlength{\rightmargin}{\leftmargin}
        \setlength{\itemsep}{-12pt}
        \setlength{\parsep}{20pt}
        \begin{lrbox}{\@tempboxa}%
        \begin{minipage}{\linewidth-2\fboxsep}
        }%
        {
        \end{minipage}%
        \end{lrbox}%
        \fcolorbox{black}{shade}{\usebox{\@tempboxa}}\newline
        }%
\makeatother

\renewcommand{\eqref}[1]{\hyperref[#1]{(\ref*{#1})}}

\newcommand{\dd}{\mathrm{d}}

\renewcommand{\dd}{{\rm d}}

\newcommand{\s}{\mathcal{S}}

\newcommand{\1}{\mathbf{1}}
\newcommand{\re}{\mathbb{R}}
\newcommand{\bP}{\mathbf{P}}

\newcommand{\bE}{\mathbf{E}}

\newcommand{\bQ}{\pmb{\rm Q}}
\newtheorem{theorem}{Theorem}
\newtheorem{lemma}{Lemma}

\newtheorem{definition}{Definition}

\newcommand*{\pref}[1]{\hyperref[#1]{(\ref*{#1})}}
\newcommand*{\refpref}[2]{\hyperref[#2]{\ref*{#1}(\ref*{#2})}}

\newcommand{\R}{\mathbb{R}}%Reales
\allowdisplaybreaks 
\newcommand{\one}[1]{\mathbbm{1}_{\{#1\}}}

\definecolor{amethyst}{rgb}{0.6, 0.4, 0.8}
\definecolor{shitbrown}{rgb}{0.55, 0.71, 0.0}
\definecolor{aqua}{rgb}{0.0, 1.0, 1.0}
\definecolor{asparagus}{rgb}{0.53, 0.66, 0.42}
\definecolor{amber(sae/ece)}{rgb}{1.0, 0.49, 0.0}
 	\definecolor{armygreen}{rgb}{0.29, 0.33, 0.13}
	\definecolor{shitbrown}{rgb}{0.43, 0.21, 0.1}
	\definecolor{brightpink}{rgb}{1.0, 0.0, 0.5}
	\definecolor{brightube}{rgb}{0.82, 0.62, 0.91}
	 	\definecolor{byzantine}{rgb}{0.74, 0.2, 0.64}
		\definecolor{chartreuse(web)}{rgb}{0.5, 1.0, 0.0}

\title{The strong law of large numbers and a functional central limit theorem for general Markov additive processes}

\author{A. E. Kyprianou\thanks{
Department of Statistics
University of Warwick
Coventry
CV4 7AL, UK. E-mail: \texttt{andreas.kyprianou@warwick.ac.uk}}
\, and V. Rivero\thanks{Centro de Investigaci\'on
 en Mathem\'aticas, A. C. 
 Calle Jalisco s/n. 36240
M\'exico.
E-mail: \texttt{rivero@cimat.mx}},    
}

\date{\it\small Dedicated to Narahari Umanath Prabhu}

\begin{document}

\maketitle
\vspace{-0.5cm}
\begin{abstract}
In this note we re-visit the fundamental question of the strong law of large numbers and central limit theorem for processes in continuous time with conditional stationary and independent increments. For convenience we refer to them as Markov additive processes, or MAPs for short. Historically used in the setting of queuing theory,  MAPs have often been written about when the underlying  modulating process is an ergodic Markov chain on a finite state space, cf. \cite{AsmussenQueue, AA} not to mention the classical contributions of Prabhu \cite{pachecoprabhu,Prabhu}. Recent works have addressed the strong law of large numbers when the underlying modulating process is a general Markov processes; cf. \cite{TAMS, Mine}. We add to the latter with a different approach based on an ergodic theorem for additive functionals and on the semi-martingale structure of the additive part. This approach also allows us to deal with the setting that the modulator of the MAP is either positive or null recurrent. The methodology additionally inspires a CLT-type result.
\medskip


\medskip

\end{abstract}

\section{Introduction} 
Suppose $((\xi_{t},\Theta_{t}), t\ge 0)$ is a bivariate coordinate process on ${\R\times\s}$, where $\s$ is a polish space and
\[
\left((\xi,\Theta),\bP \right)=\left((\xi_{t},\Theta_{t}, t\ge 0),\mathcal{F}_{\infty},(\mathcal{F}_{t}, t\ge 0),\{\bP_{x,\theta}:(x,\theta)\in\R\times\s\}\right),
\]
 is a (possibly killed) Markov process with $\bP_{x,\theta}\left(\xi_{0}=x,\Theta_{0}
=\theta\right)=1$. Here $(\mathcal{F}_{t}, t\ge 0)$ is the minimal augmented and right continuous admissible filtration, and $\mathcal{F}_{\infty}=\bigvee_{t=0}^{+\infty}\mathcal{F}_{t}$.

  \begin{definition}
    The process $\left((\xi,\Theta),\bP\right)$ is called a Markov additive process (MAP) on $\R\times\s$ if, for any $t \geq 0$, given $\{(\xi_s,\Theta_s), s\leq t\}$, the process $(\xi_{s+t}-\xi_t,\Theta_{s+t})_{s\geq 0}$ has the same law as $(\xi_{s},\Theta_{s}, s\ge 0)$ under $\bP_{0,v}$ with $v=\Theta_t$. We call $\left((\xi,\Theta),\bP\right)$ a nondecreasing MAP if $\xi$ is a nondecreasing process on $\R$.
  \end{definition}

For a MAP  $\left((\xi,\Theta),\bP\right)$, we call $\xi$ the \textit{ordinate} and $\Theta$ the \textit{modulator}. By definition we can see that a MAP is conditionally translation invariant in $\xi$ in the sense that $((\xi_{t},\Theta_{t}, t\ge 0),\bP_{x,\theta})$ is equal in law to $((\xi_{t}+x,\Theta_{t}, t\ge 0),\bP_{0,\theta})$ for all $x\in \R$ and $\theta\in\s $. In words: the `conditional' translation is dependent on the current value of the modulator.  We assume throughout the paper the following path regularity assumptions for the process $(\xi,\Theta)$.

\begin{itemize}
\item[{\bf (H1)}] The process $((\xi_t, \Theta_{t}), t\ge 0)$ is a $\R\times\s$-valued Hunt process. 
\end{itemize}

 Whilst MAPs have found a prominent role  in e.g. classical applied probability models for queues and dams, c.f. \cite{Prabhu, AsmussenQueue} when $\Theta$ is a Markov chain, the case that $\Theta$ is a general Markov process has received somewhat less attention. Nonetheless, a core base of literature exists in the general setting from before the year 2000 thanks to e.g. \cite{arjasspeed,AsmussenQueue,AA,CPR,Cinlar1, CinlarI&II, CinlarMRT,KPal, Kaspi1,Kaspi2,NeyNummelin,NNC,NNSLT} and in  the recent century due to their intimate relationship with self-similar Markov processes and multi-type scaling limits, see e.g. \cite{LCGZ,KKPW,Stephenson,SilvaPardo} and the reference therein.

\medskip

Fundamentally the ordinate of a MAP has similar qualitative long-time behaviour to L\'evy processes and random walks, and extends that of additive functionals of Markov processes, the most simple of which is the strong law of large numbers. That is, the almost sure existence of the limit $\lim_{t\to\infty}\xi_t/t$. Until recently results concerning the strong law of large numbers for MAPs were largely restricted to the setting that $\s$ is a finite set; cf. \cite{AsmussenQueue, AA}. 

\medskip

In addition to (H1), a reasonable assumption in order to derive a law of large numbers is that the underlying modulator process $\Theta$ has a recurrent structure. A natural setting is therefore captured with our second assumption.

\begin{itemize}
\item[{\bf (H2)}] The process $(\Theta_{t}, t\ge 0)$ is Harris recurrent with a $\sigma$-finite invariant measure $\pi.$ We will say that the process is positive recurrent if $\pi$ is a finite measure, and hence can be normalized to be a probability measure; and we will say it is zero recurrent if $\pi$ is an infinite  measure.
\end{itemize}

\noindent Sufficient conditions in terms of $(\xi,\Theta)$ for $\pi$ to exist are given in \cite{TAMS}. There, the set $\mathcal{S}$ is taken as a subset of $\mathbb{S}^{d-1}$. Moreover, in place of (H2), it was assumed that there exist a constant $\delta>0$, a probability measure $\rho$ on $\mathcal{S}$ and  a family of measures $\{\phi(\theta,\cdot):\theta\in\s \}$ on $\R$ with $\inf_{\theta\in\mathcal{S}}\phi(\theta,\R)>0$ such that
  \begin{equation}\label{condi:Harris type}
    \bP_{0,\theta}\left(\xi_{\delta}\in \Gamma,\, \Theta_{\delta}\in A\right)\ge \phi(\theta,\Gamma)\rho(A)\quad\forall\theta\in\mathcal{S},\ A\in\mathcal{B}(\mathcal{S}),\ \Gamma\in\mathcal{B}(\R).
    \end{equation}
    which is a sufficient condition to assume (H2) holds.

\medskip
Together with a relatively strong condition on first moments of $\xi$, \cite{TAMS} established a strong law of large numbers (SLLN). More precisely,   remembering that $\mathcal{S}$ was chosen as a subset of $\mathbb{S}^{d-1}$,  under the assumptions $\bE_{0,\pi}[\sup_{t\in[0,1]}|\xi_t|]<\infty$ and \eqref{condi:Harris type}, it was shown in Proposition 2.15  of \cite{TAMS} that 
\[
\lim_{t\to\infty}\frac{\xi_t}{t} = \bE_{0,\pi}[\xi_1].
\]
The proof of this SLLN in \cite{TAMS} is based around the representation of the invariant measure $\pi$ as an occupation measure over regeneration cycles.
 Related results concerning the long term behaviour of the limsup and liminf of $\xi$ were recently shown in \cite{Mine} and \cite{DDK}.
 
 \medskip

 In this article we prove a version of the SLLN under conditions which are based on the characteristics of the process as well as accommodating for the setting that $\pi$ has infinite mass (the case that $\Theta$ is null-recurrent). Our approach uses known results concerning the ergodic behaviour of additive functionals for Harris recurrent Markov processes as well as the semi-martingale representation of the ordinator process in combination with  known SLLN-type results for semi-martingales. 
 \medskip

  The method we used for the SLLN inspires further results concerning functional Central Limit Theorem (CLT) type results, again the same two ingredients being key. The nature of the CLT we will present here for MAPs will be to find an additive process $A$ and a function $h$ such that 
  \[
\frac{1}{\sqrt{h(n)}}\left(\xi_{nt}-A_{nt}\right)
  \]
  converges on the Skorokhod space as $n\to\infty$ to an anomalous diffusion. 
  
  Earlier results of this type have arisen in the queueing literature,  typically under the assumption that the modulator is a positive recurrent Markov chain, \cite{Mandjesetal,Pang,Ankita,steichen}. Here we allow the modulator to be any Harris recurrent Markov processes, null or positive recurrent. This level of generality is achieved by exploiting existing results from the vast literature on functional CLT for semi-martingales and additive functionals of Markov processes. 
  \medskip

 In order to state our main results, we need to first recall some facts concerning the semimartingale structure of MAPs as laid out in \cite{CinlarMRT, CinlarI&II, Cinlar1}.
\section{MAP as semi-martingales}
In a series of seminal papers, \c{C}inlar laid the foundations for using  stochastic calculus to describe the ordinator of a MAP as a functional of the modulator.   We will summarise some of his results below. By identifying MAPs as processes with conditional independent increments (C-PII), this will also allow us to extract some results from Section 2.6 in \cite{J&S}, where general C-PII processes are studied.  
\medskip

We denote by $\mathcal{G}=(\mathcal{G}_t, t\geq 0)$  the natural $\sigma$-field and filtration, respectively, generated by $\Theta$,   completed with respect to the family  $\int_{\theta\in\s}\mu(\dd\theta)\bP_{0,\theta}(\cdot)$ for $\mu$ a (probability) measure on $\s.$ Proposition 2.20 in part II of \cite{CinlarI&II} establishes the existence of a regular version of the conditional law of $\xi$ given $\mathcal{G}$. We will denote it by $\left(\bQ^{\omega}, \omega\in\mathbb{D}([0,\infty), \s)\right);$ where as usual  $\mathbb{D}\left([0,\infty), \s\right)$ denotes the set of $\s$-valued cadlag paths with a lifetime. For any bounded paths functionals $F$ and $H$ on  $\mathbb{D}([0,\infty), \s)$ and $\mathbb{D}([0,\infty), \mathbb{R})$, respectively, with the obvious meaning for $\mathbb{D}([0,\infty), \mathbb{R})$, we have 
\begin{equation*}
\begin{split}
&\bE_{0,\theta}\left[{F}(\Theta_t, t\geq 0) {H}(\xi_t, t\geq 0)\right]\\
&=\bE_{0,\theta}\left[{F}(\Theta_t, t\geq 0) \bQ^{\omega}\left[{H}(\xi_t, t\geq 0)\right]\Large|_{\omega=(\Theta_t, t\geq 0)}\right],\qquad \forall\theta\in\s.
\end{split}
\end{equation*}
When conditioning with respect to $\mathcal{G}$ we will systematically work with the aforementioned regular version. 
\medskip

Theorem 2.22 in part II of \cite{CinlarI&II} establishes that for any given $\omega\in\mathbb{D}([0,\infty), \s),$  $\xi$ has independent increments under $\bQ^{\omega}$.  Otherwise said, conditionally on $\Theta,$ $\xi$ has independent increments. Theorem 2.23 of the same reference, proves that there exist stochastic processes $\xi^{f}$, $\xi^{c}$ and $\xi^{d}$, which are conditionally independent given $\mathcal{G}$, and an additive functional of $\Theta,$ say $\chi,$ which does not  necessarily have finite variation over finite intervals, such that
 \begin{equation}
 \xi_t=\chi_t+\xi^{c}_t+\xi^{f}_t+\xi^{d}_t,\qquad t\geq 0.
 \label{Cinlardecomp}
 \end{equation}
Moreover,  these processes are such that:
\begin{itemize}
\item $\xi^{c}=(\xi^{c}_t, t\geq 0)$ is a continuous process, $( \xi^{c}, \Theta)$ is a MAP and conditionally on $\mathcal{G},$ $\xi^{c}$ is a Gaussian process; 
\item $\xi^{f}=(\xi^{f}_t, t\geq 0)$ is a pure jump process, $( \xi^{f}, \Theta)$ is a MAP, and if $\tau$ is jump time of $\xi^{f},$ then $\tau$ is $\mathcal{G}$-measurable, that is, the discontinuities of $\xi^{f}$ are determined by $\Theta;$
\item $\xi^{d}=(\xi^{d}_t, t\geq 0),$ is such that $( \xi^{d}, \Theta)$ is a MAP, and, conditionally on $\mathcal{G}$, $\xi^d$  is a stochastically continuous process with independent increments.
\end{itemize}
This decomposition is reminiscent of the L\'evy-It\^{o} description of additive processes, see e.g. \cite{Satobook} Section 19. We can think of $\xi^{d}$ as an additive Cox process, i.e. conditionally on $\Theta,$ this process is constructed as in  the proof of Theorem 19.2 in~\cite{Satobook} using a time inhomogeneous Poisson process whose characteristics depend on $\Theta$. The process $\chi$ plays the role of a drift, but its paths are allowed to be rather irregular. Analogously, $\xi^{c}$ is morally a stochastic integral of a predictable process of $\Theta$ with respect to a Brownian motion. 

In addition to this, we also note the following observations concerning other decompositions that have been noted in the literature. On the one hand, from the conditional version of Theorem 5.1 in Chapter II in \cite{J&S}, page 114, we know that $\xi$ can be written as 
$\widetilde{\xi}+\xi^{\prime}$ where $\xi^{\prime}$ is a process with conditional independent increments and is a semi-martingale; and $(\widetilde{\xi}_{t}, t\geq 0),$ is a c\'adl\'ag process issued from zero and adapted to the filtration of $\Theta.$  On the other hand, in addition to the decomposition \eqref{Cinlardecomp}, Theorem 6.5 in Chapter II in \cite{J&S}, shows that $\xi$ is a semi-martingale if and only if for each $\lambda\in\R$ the process
$$t\mapsto\bE\left[e^{i\lambda \xi_t} \large| \mathcal{G}\right], \qquad t\geq0,$$ admits a version that is c{\`{a}}dl{\`{a}}g with finite variation over bounded intervals. Corollary 5.11 in Chapter II of \cite{J&S} also implies that, conditionally on $\Theta,$ $\xi$ is a semimartingale if and only if $\chi+\xi^{f}$ is a process with finite variation over finite intervals. 
\medskip

We thus conclude these observations by arguing that there is  no major  loss of generality in assuming that 
\begin{itemize}
\item[{\bf (H3)}] $\xi$ is a semimartingale.
\end{itemize}

In \cite{CinlarI&II}, \c Cinlar also gave a description of the characteristic function of $\xi$ conditionally on $\mathcal{G},$ and hence described the characteristics of  $\xi$ as a semimartingale. To recall that description, we quote first another result due to \c Cinlar that describes the jump structure of $(\Theta, \xi).$  Indeed, in \cite{Cinlar1} it is proved that  there exists a continuous increasing additive functional $t\mapsto H_{t}$ of $\Theta$ and a transition kernel $\Pi$ from $\s $ to $\s \times \R$ satisfying that, 
for every nonnegative predictable process $Z$ and nonnegative measurable function $g:\s \times \R \times \s \times \R\to \R^{+}$, we have the {\it compensation formula} 
\begin{align}
\bE_{0,\theta}&\left[\sum_{s\le t}Z_{s}g(\Theta_{s-},\xi_{s-},\Theta_{s},\xi_{s})\1_{\{\Theta_{s-}\neq \Theta_{s}\ \text{or}\ \xi_{s-}\neq \xi_{s}\}}\right]\nonumber\\
&=\bE_{0,\theta}\left[\int^{t}_{0}{\rm d}H_{s}Z_{s}\int_{\s \times \R}\Pi(\Theta_{s}, {\rm d}v, {\rm d}y)g(\Theta_{s},\xi_{s},v,\xi_{s}+y)\right]\label{levysys}
\end{align}
for all $\theta\in\s $ and $t\ge 0$.  The pair $(H,\Pi)$ is said to be a \textit{L\'{e}vy system} for $\left((\xi,\Theta),\bP\right)$. 
\begin{itemize}
\item[{\bf (H4)}]We will assume that $H_t=t\wedge{\zeta},$ $t\geq0.$ 
\end{itemize}
This assumption is indeed satisfied in many explicit examples that appear in the literature. When it is  not satisfied, (H4) can be fulfilled by applying  the time change $t\mapsto H^{-1}_t$ to the original MAP. As such, (H4) is not a major restriction. 
\medskip

As mentioned before, there are some jumps of $\xi$ that are induced by $\Theta,$ which means that when we condition on the path of $\Theta,$ some fixed discontinuities appear in the path of $\xi,$ through $\xi^{f}$ and  $\chi.$   
\medskip

Furthermore, \c Cinlar proved the following properties on the jump kernel $\Pi$ 
\begin{itemize}
\item[(i)] $\Pi(\theta, \{\theta\}, \{0\})=0,$ for all $\theta\in\s;$
\item[(ii)] $\int_{\R\setminus\{0\}}\Pi(\theta, \{\theta\}, \dd y)(1\wedge|y|^{2})<\infty,$ for all $\theta\in \s;$ 
\item[(iii)] there exists a kernel $K,$ from $\s$ into $\s,$ and a sub-probability kernel $\nu,$ from $\s^{2}$ into $\R,$ such that 
$$\Pi(\theta, \dd \beta, \dd y)\one{\beta\neq \theta}=K(\theta, \dd\beta)\one{\theta\neq\beta}\nu_{\theta,\beta}(\dd y),\theta,\beta\in \s, \ y\in\R;$$ moreover, $(H, K)$ is a L\'evy system for $\Theta.$
\end{itemize}
For notational convenience, we define a random measure on $[0,\infty)\times\R\setminus\{0\}$ by 
\begin{equation}\label{L}
L(\dd s, \dd x)= \Pi(\Theta_s, \{\Theta_s\}, \dd x )\dd s +\sum_{u>0}\one{\Theta_u\neq \Theta_{u-}}\delta_{u}(\dd s)\nu_{(\Theta_{u-}, \Theta_{u-})}
(\dd x).\end{equation}

Corollary 2.25 in part II of \cite{CinlarI&II}, Lemma 2.24 in \cite{Cinlar1} and an application of the conditional version of the identity (4.16) and Theorem 5.2,  both from Chapter II in~\cite{J&S} (see also Theorem 6.6 in Chapter II in \cite{J&S} for the semi-martingale case), imply that  
\begin{equation}\label{eq:charac}  
\begin{split}
&\bQ^{\omega}\left[\exp\{i\lambda(\xi_{t}-\xi_{s})\}\right]\\
&=\exp\left \{i\lambda(B_t(\omega)-B_s(\omega))-\frac{\lambda^{2}}{2}(C_{t}(\omega)-C_{s}(\omega))\right.\\ &\quad +\left.\int^{t}_{s}\dd u\int_{\R\setminus\{0\}}\left(e^{i\lambda x}-1-\lambda x\one{|x|<1}\right)\Pi(\omega(u), \{\omega(u)\},\dd x)\}\right\}\\
&\quad\times\prod_{s<u\leq t}e^{-i\lambda \Delta B_u(\omega)}\left(1+\one{\omega(u-)\neq \omega(u)}\int_{\R}(e^{i\lambda x}-1)\nu_{(\omega(u-), \omega(u))}(\dd x)\right);
\end{split}
\end{equation}
and 
\begin{equation*}
\begin{split}
\bQ^{\omega}\left[\Delta\xi_{t}\in\dd y\right]=\nu_{(\omega(t-), \omega(t))}(\dd y)\one{y\neq 0}+\left(1-\nu_{(\omega(t-), \omega(t))}(\R)\right)\delta_{0}(\dd y),\qquad t\geq0;
\end{split}
\end{equation*}
where $C=(C_t, t\geq 0)$ is a positive valued continuous additive functional of $\Theta.$ This corresponds to the quadratic variation of $\xi^{c}$ given $\mathcal{G}$, so it is such that $$(\xi^{c}_t)^{2}-C_t,\qquad t\geq 0,$$ is a local martingale under $\bQ^{\cdot}.$ Furthermore, $B$ is the compensator of  $\xi^{1}=(\xi^{1}_t, t\geq 0),$ the process which is obtained from $\xi$ by removing the jumps of magnitude larger than unity. 
\medskip

In order to prove the SLLN, we will require that for each $t>0,$ $\bQ^{\omega}[|\xi_t|]<\infty,$ which happens if and only if 
\begin{equation}\label{M1}
\begin{split}
&\int^{t}_{0}\dd u\int_{\R\setminus\{0\}}|x|\one{|x|>1}\Pi(\omega(u), \{\omega(u)\},\dd x)\\
&\qquad+\sum_{0<s\leq t}\int_{\R}|y|\nu_{(\omega(s-), \omega(s))}(\dd y)\one{y\neq 0}<\infty. 
\end{split}
\end{equation}
In what follows, we will assume that \eqref{M1} is satisfied, however, as we will see,  in our main result, it will be implied by a stronger assumption, to be introduced in (H7). In which case, 
\begin{equation}\label{compmean}
\begin{split}
\bQ^{\omega}[\xi_t]&=B_t-\sum_{0<s\leq t}\Delta B_s+\int^{t}_{0}\dd u\int_{\R\setminus\{0\}}x\one{|x|>1}\Pi(\omega(u), \{\omega(u)\},\dd x)\\
&\qquad+\sum_{0<s\leq t}\int_{\R}y\nu_{(\omega(s-), \omega(s))}(\dd y)\one{y\neq 0}<\infty.
\end{split}
\end{equation}
Additionally, when \eqref{M1} holds, Proposition 2.29 in Chapter II in~\cite{J&S} implies that  $\xi$ is a semimartingale under $\bQ^{\omega}$ and its canonical decomposition takes the form  
\begin{equation}\label{xisemi-mart}
\xi_t=\xi_0+N_t+\widetilde{B}_t, \qquad t\geq0,
\end{equation}
 where $N$ is a locally square integrable martingale and $\widetilde{B}$ is such that 
\begin{equation}\label{Btilde}
\begin{split}
\bQ^{\omega}[\xi_t]=\widetilde{B}_t=B^{c}_{t}+&\int^{t}_{0}\dd u\int_{\R\setminus\{0\}}x\one{|x|>1}\Pi(\omega(u), \{\omega(u)\},\dd x)\\
&+\sum_{0<s\leq t}\int_{\R}y\nu_{(\omega(s-), \omega(s))}(\dd y)\one{y\neq 0},
\end{split}
\end{equation}
 with $B^{c}$ the process defined by
\begin{equation}
B^{c}_t:=B_t-\sum_{0<s\leq t}\Delta B_s,\qquad t\geq 0.
\end{equation}
Proposition 2.9 in Chapter II in \cite{J&S} allow us to ensure that the processes $B^{c},$ and $C$ are absolutely continuous with respect to $H.$ Hence, there exists predictable processes $b$ and $c$ such that $$B^{c}_t=\int^{t}_{0}b_s \dd s,\qquad C_t=\int^{t}_{0}c_s \dd s,\qquad t\geq 0.$$ 
Taking expectations in the left hand side of \eqref{M1} and using that $(H,K)$ is a L\'evy system for $\Theta,$ we have that 
\begin{equation}\label{display9}
\begin{split}
&\bE_{0,\theta}\left[\int^{t}_{0}\dd u\int_{\R\setminus\{0\}}|x|\one{|x|>1}\Pi(\Theta_u, \{\Theta_u\},\dd x)+\sum_{0<s\leq t}\int_{\R}|y|\nu_{(\Theta(s-), \Theta(s))}(\dd y)\one{y\neq 0}\right]\\
&=\bE_{0,\theta}\left[\int^{t}_{0}\dd u\int_{\R\setminus\{0\}}|x|\one{|x|>1}\Pi(\Theta_u, \{\Theta_u\},\dd x)\right.\\
&\qquad\hspace{1cm} \left. +\int^{t}_{0}\dd s \int_{\s}K(\Theta_{s}, \dd\beta)\int_{\R}|y|\nu_{(\Theta_{s}, \beta)}(\dd y)\one{y\neq 0}\right]\\
&=\bE_{0,\theta}\left[\int^{t}_{0}\dd u\int_{\R\setminus\{0\}}|x|\one{|x|>1}\Pi(\Theta_u, \s,\dd x)\right.\\ &\qquad\hspace{1cm}+\left.\int^{t}_{0}\dd u\int_{\R\setminus\{0\}}|x|\one{|x|\leq 1}\Pi(\Theta_u, \s\setminus\{\Theta_u\},\dd x)\right].
\end{split}
\end{equation}
We will use below that $\xi$ has a finite mean, this will be implied by the assumption
\begin{itemize}
\item[\bf (H5)] For every $t>0,$ 
\begin{equation*}
\begin{split}
&\bE_{0,\theta}\left[\int^{t}_{0}|b_s| \dd s + \int^{t}_{0}\dd u\int_{\R\setminus\{0\}}|x|\one{|x|>1}\Pi(\Theta_u, \s,\dd x)
\right.\\ &\qquad\qquad+\left.\int^{t}_{0}\dd u\int_{\R\setminus\{0\}}|x|\one{|x|\leq 1}\Pi(\Theta_u, \s\setminus\{\Theta_u\},\dd x)\right]
<\infty.
\end{split}
\end{equation*}
\end{itemize}
In that case, we will use the notation:
\begin{equation}\label{M2}
\mu^{d}(\theta):=\int_{\R\setminus\{0\}}x\one{|x|>1}\Pi(\theta, \s,\dd x)+\int_{\R\setminus\{0\}}x\one{|x|\leq 1}\Pi(\theta, \s\setminus\{\theta\},\dd x),\qquad \theta\in\s.
\end{equation}
We will denote by $A$ the additive functional defined by 
\begin{equation}\label{Adef}
A_t:=B^{c}_t+\int^{t}_0\mu^{d}(\Theta_s)\dd s,\qquad t\geq 0.
\end{equation}
Notice that arguing as in \eqref{display9} it is proved that process $A$ is the predictable compensator of the process $(\widetilde{B}_{t}=\bQ^{\Theta}[\xi_t], t\geq 0),$ which is described in \eqref{Btilde}.
\medskip

To state our next assumption, we recall that the Revuz measure $\nu_{\alpha}$ of an additive functional $\alpha=(\alpha_t, t\geq 0),$ of $\Theta,$ is defined via the formula
$$\nu_{\alpha}(f)=\lim_{t\to0}\frac{1}{t}\bE_{0,\pi}\left[\int^{t}_{0}f(\Theta_s)\dd \alpha_s\right],$$
 for $f:\s\mapsto\re^{+}$. We will denote its total variation by $||\nu_{\alpha}||.$ See e.g. \cite{RY} Section X.2 for background.

\begin{itemize}
\item[\bf (H6)] The total variation $||\nu_{C}||$ of the Revuz measure of $C$ is finite, and 
$$\int_{\s}\pi(\dd\theta)\int_{|x|\leq 1}\Pi(\theta, \{\theta\},\dd x)x^{2}<\infty.$$
\item[\bf (H7)] The total variation $||\nu_{B^{c}}||$ of the Revuz measure of $B^{c}$ is finite and  
 $$\int_\s \pi(\dd \theta)\mu^{d}(\theta)<\infty.$$ In that case, we will write
 $$m_1:=||\nu_{B^{c}}||+\int_\s \pi(\dd \theta)\mu^{d}(\theta).$$
 \item[\bf (H8)] There is a $p>0$ such that
 $$\int_{\s}\pi(\dd\theta)\int_{x\in\R}\Pi(\theta, \s, \dd x)|x|^{1+p}\one{|x|>1}<\infty.$$
\end{itemize}
Observe that (H7) implies (H5), and that similar computations to the ones above (H5) allow us to prove, using \eqref{compmean}, that 
$$m_1=\bE_{0,\pi}[\xi_1].$$

Finally we note that 
under the assumptions (H1-5) we have that
\[
\tilde{\xi}_t: = \xi_t -A_t, \qquad t \geq0,
\]
is  a centered MAP; that is, $(\tilde\xi, \Theta)$ is a MAP with zero mean $$\bE_{x,\theta}[\tilde\xi_t] = 0,\qquad t\geq0, \qquad x\in\R,\theta\in\s.$$ 
More so, appealing to the property of conditional stationary and independent  increments, 
\[
\bE_{x,\theta}[\tilde\xi_{t+s}| \mathcal{F}_{t}] = \tilde\xi_t + \bE_{x',\theta'}[\tilde\xi_s]_{x' = \xi_t, \theta'=\Theta_t} = \tilde\xi_t, \qquad s,t\geq0.
\]
This  tells us that although $(\tilde\xi_t, t\geq0)$ is not a Markov process, it is a martingale, which is the first claim of the following  Lemma; the remaining claims are easy to verify from  
\eqref{eq:charac}.

\begin{lemma}\label{square} 
Under the assumptions (H1-5), $\xi-A$ is a martingale. If moreover (H8) is satisfied with $p=1,$ we have that $\xi-A$ is a $\mathcal{F}$-square integrable martingale under $\bP$, and its quadratic variation is the additive functional of $\Theta$
$$<\xi-A>_t=C_t+\int^{t}_{0}\dd s\int_{\R\setminus\{0\}}x^{2}\Pi(\Theta_s, \s, \dd x),\qquad t\geq0.$$  
\end{lemma}

\section{Main results}
Having discussed  various aspects of the semimartingale representation of $(\xi, \Theta)$, we can now state our main results.

\begin{theorem}[Strong law of large numbers for MAP]\label{SLLNMAP}
Suppose $\pi$ is an infinite measure, that is, $\Theta$ is null-recurrent, then under the assumptions (H1-4), and (H6-8), we have that 
\begin{equation}
\frac{\xi_t}{A_t}\xrightarrow[t\to\infty]{\text{a.s.}} 1 \text{ and } \frac{A_t}{t}\xrightarrow[t\to\infty]{\text{a.s.}} 0
\label{tildaSLLN+}
\end{equation}
if $m_1 \neq 0$ and otherwise if $m_1 = 0$ then 
\begin{equation}
\frac{\xi_t}{t}\xrightarrow[t\to\infty]{\text{a.s.}} 0  \text{ and } \frac{A_t}{t}\xrightarrow[t\to\infty]{\text{a.s.}} 0.
\label{tildaSLLN++}
\end{equation}
When $\pi$ is a probability measure, that is when $\Theta$ is positive recurrent, then under the assumptions (H1-4) and (H7), we have that
\begin{equation}
\frac{\xi_t}{t}\xrightarrow[t\to\infty]{\text{a.s.}} \bE_{0,\pi}[\xi_1] \text{ and } \frac{A_t}{t}\xrightarrow[t\to\infty]{\text{a.s.}} \bE_{0,\pi}[\xi_1].
\label{tildaSLLN*}
\end{equation}
\end{theorem}

There are a number of remarks we can make about the above SLLN. First, the main contribution is the setting that $\pi$ is an infinite measure. Indeed, as we will see from the proof, it is relatively straightforward to deduce the second part of the theorem from what is already known in Proposition 2.15 of \cite{TAMS}.

\medskip

Second, the assumption (H7) is the natural extension of the usual assumption for L\'evy processes for the SLLN to hold. In the L\'evy setting,  the assumption (H6) is  automatically satisfied. In the setting of MAPs some control of the volatility and linear terms is required in the long run, and this involves the  stationary regime of the modulator.  In essence, the conditions ensure that  there are no states of the modulator that would lead to an extreme variation in the ordinator, causing the SLLN not to hold.   
\medskip

%With the semi-martingale method we can prove the SLLN in the case where the process is Harris's null recurrent, that means that the measure $\pi$ has infinite mass. In that case, the normalization is $\xi_t/A_t$. 

 Following the SLLN it is natural to ask whether a central limit theorem holds. To answer this question positively, we need to recall the definition of a Mittag-Leffler process. 
\medskip

For $\alpha\in(0,1),$ let $\sigma^{(\alpha)}$ be an $\alpha$-stable subordinator, normalized so to have $$\log\bE[e^{-\sigma^{(\alpha)}_1}]=1.$$ For $\alpha=1$ we take the pure drift process $\sigma^{(1)}_{t}=t,$ $t\geq 0.$ For any $\alpha\in(0,1],$ we define the process of the first passage time of $\sigma^{(\alpha)}$ by taking its right-continuous inverse
$$W^{(\alpha)}_t=\inf\{s>0: \sigma^{(\alpha)}_s>t\},\qquad t\geq 0.$$ We will refer to this process as a Mittage-Leffler process of index $\alpha\in(0,1].$ 
We now have all the elements required to state our CLT.

\begin{theorem}\label{CLTMAPS}
Suppose that  (H1-4) and (H6-8) hold with $p=1.$  Moreover, we  assume the so-called  Lindenberg condition that, for all $\epsilon>0$,
 \begin{equation}\label{lin}
\frac{1}{t}\int_{(0,t]\times\R\setminus\{0\}}L(\dd u, \dd x)x^{2}\one{|x|>\epsilon\sqrt{t}}\xrightarrow[t\to\infty]{\text{a.s.}}0.
\end{equation}  
In addition, we assume that
 there is a regularly varying function  $h:(0,\infty)\mapsto(0,\infty)$ at infinity, with some index $\alpha\in(0,1],$ such that for every measurable function  $g:\s\mapsto\re^{+}$ such that $0<\pi(g)<\infty,$ 
$$\bE_{0,\theta}\left[\int^{\infty}_{0}e^{-s/t}g(\Theta_s)\dd s\right]\sim h(t)\pi(g),\qquad t\to\infty;$$
for $\pi$-almost all $\theta\in\s;$
\medskip

With all these assumptions in place, we have  that, under $\bP_{0,\theta},$ there is the following convergence in the sense of Skorohod's topology
$$\left(\frac{1}{\sqrt{h(n)}}\left(\xi_{nt}-A_{nt}\right), \ t\geq 0\right)\xrightarrow[n\to\infty]{}\left(\left(\nu_{\langle\xi\rangle}+2J\right)^{1/2}\Sigma_{W^{(\alpha)}_t}, \ t\geq 0\right);$$ where $A$ is defined in \eqref{Adef}, $\Sigma_{W^{(\alpha)}}$ is the anomalous diffusion obtained by the subordination of $\Sigma,$ a standard Brownian motion, subordinated by an independent Mittag-Leffler process of parameter $\alpha,$ $W^{(\alpha)},$ and \begin{equation}\label{J}
J:=\int_{\s}\pi(\dd \theta)\int_{\s}K(\theta, \dd \beta)\left(\int_{\R}\nu_{\theta, \beta}(\dd x)x\right)^{2}<\infty.
\end{equation} 
\end{theorem}
 Lindenberg's condition is commonly required when establishing CLT outside the i.i.d. case. Fundamentally, MAPs support  two main sources of randomness, that is $\xi$ and $\Theta$. To control either, one needs some moment conditions, such as (H8), and the recurrence structure of $\Theta.$ We have allowed $\Theta$ to be either positive recurrent or null recurrent. This  determines the time it takes for $\Theta$  to reach stability and complete a recurrence epoch. Such recurrent epochs are intimately connected to the increments of $\xi$ and we need to have a good understanding of them. The condition on the regular variation of the resolvent of $\Theta$ is a version of the Darling-Kac condition, introduced in \cite{HandL} extending \cite{Darling-kac}, that allows one to ensure that the sum of recurrence epochs is in the domain of attraction of a stable subordinator. The Mittag-Leffler clock appears from estimating the rate of increase of the number of recurrence epochs completed as time goes to infinity. Different behaviour arises in the positive recurrent case and the null recurrent case. 

When the process is positive recurrent, $\pi$ is finite, the assumptions in the latter Theorem are satisfied with $h(t)=t,$ for all $t\geq 0,$ and hence $\alpha=1$, see e.g. \cite{HandL} Theorem 3.1 and the discussion following it. This can be easily seen by decomposing into excursions, but for brevity we omit the details. Thus, the assumption in the previous theorem is only required when the process is null-recurrent, in that case, one might have or not the existence of the norming function $h,$ in general it will be related to the intensity rate at which long recurrence epochs appear. For background on the Darling-Kac theory we refer to~\cite{Darling-kac,Athreya,BGT,HandL}.     

The proof of Theorem~\ref{CLTMAPS} relies on the additive structure of $\xi$ using semi-martingale techniques. In the proof, the contribution from the different sources of randomness will be made apparent, and this will explain the appearance of the constants $\nu_{\langle\xi\rangle}$ and $J$. 

The remainder of the paper is devoted to prove our two main results.

  \section{Proof of Theorem \ref{SLLNMAP}}
We will assume that $m_1: = ||\nu_{B^{c}}||+\int_\s \pi(\dd \theta)\mu^{d}(\theta)\neq 0,$ and without loss of generality that $m_1>0.$ Later we will explain how to modify the proof when $m_1=0.$ Let us assume that 
$\pi$ is an infinite  measure. 
\medskip

The Ergodic Theorem of X.3.12 of \cite{RY}, in the positive and null recurrent cases, implies also that if $m_1>0,$ then
\begin{equation}\label{C/A}
\frac{C_t+\int^{t}_0\dd u \int_{|x|\leq 1}\Pi(\Theta_u, \s, \dd x)x^{2}}{A_t}\xrightarrow[t\to\infty]{\text{a.s.}}\frac{||\nu_{C}||+\int_{s}\pi(\dd\theta)\int_{|x|\leq 1}\Pi(\theta, \s, \dd x)x^{2}}{||\nu_{B^{c}}||+\int_\s \pi(\dd \theta)\mu^{d}(\theta)}.
\end{equation}
The latter is finite whenever (H6) holds. 
\medskip

Recall that $
\tilde{\xi}_t: = \xi_t -A_t$, $t\geq0$,
where $A$ was defined in \eqref{Adef}.
Next, the main result in \cite{Lipster} will allow us to prove that 
\begin{equation}\label{tildaSLLNA}
\frac{\widetilde{\xi}_{t}}{A_t}\xrightarrow[t\to\infty]{\text{a.s.}}0.
\end{equation}
To  that end, we need to introduce the quantity
\[
M_t : = \int_0^t \frac{1}{1+A_s}\dd  \tilde\xi_s, \qquad t\geq0.
\]
We notice that this is a local martingale, whose quadratic variation, $\langle M\rangle_t$, $t\geq0$, satisfies 
\[
\langle M\rangle_t = \int_0^t \frac{1}{(1+A_s)^2}\dd \langle\tilde\xi\rangle_s.
\]
In particular, appealing to the semimartingale decomposition of $\xi$,  the continuous part of $\langle M\rangle$,  denoted by $\langle M^c\rangle$, satisfies
\[
\langle M^c\rangle_t = \int_0^t  \frac{1}{(1+A_s)^2}\dd C_s, \qquad t\geq0.
\]
Finally, we define the quantity
\[
G_t = \langle M^c\rangle_t + \sum_{s\leq t} \frac{(\Delta M_s)^2}{1+|\Delta M_s|} %= \langle m^c\rangle_t + \sum_{s\leq t} \frac{(\Delta \xi_s)^2}{(1+A_s)^2}\frac{1+A_s}{1+A_s+|\Delta \xi_s|} 
, \qquad t\geq0.
\]

Theorem 1 of  \cite{Lipster} now tells us that if the predictable compensator of $(G_t, t\geq0)$, written $(\widetilde{G}_t, t\geq0),$ is such that $ \widetilde{G}_\infty: = \lim_{t\to\infty}\widetilde{G}_t<\infty$ almost surely, then \eqref{tildaSLLNA} holds. 

\medskip

In order to verify this we use that the point process of jumps of $m$, $$(\Delta M_t, t>0)$$ coincides with $$\left(\frac{1}{(1+A_t)}\Delta\xi_{t}, t>0\right),$$ see e.g. page 329 in \cite{D&M}. This implies that
\[
G_t = \langle M^c\rangle_t + \sum_{s\leq t} \frac{(\Delta \xi_s)^2}{1+A_s}\frac{1}{(1+A_s+|\Delta \xi_s|)} 
, \qquad t\geq0,
\]
which tells us that its predictable compensator is given by
\begin{align}
\widetilde{G}_t &= \int_0^t  \frac{1}{(1+A_s)^2}\dd C_s + \int_0^t \dd s
\int_\s\Pi(\Theta_s, \s, \dd x)\frac{x^2}{1+A_s}\frac{1}{1+A_s + |x|}.
\label{splittheintegral}
\end{align}
Now, to study these integrals we split the second  integral on the right-hand side of \eqref{splittheintegral} according to whether the jumps are  smaller than or greater than unity in magnitude. 
We begin with the former, for that we use the convergence in \eqref{C/A} to ensure that for any $\epsilon >0,$ there is a large $T$ such that 
$$\left|\frac{C_t+\int^{t}_0\dd u \int_{|x|\leq 1}\Pi(\Theta_u, \s, \dd x)x^{2}}{A_t}-\frac{||\nu_{C}||+\int_{s}\pi(\dd\theta)\int_{|x|\leq 1}\Pi(\theta, \s, \dd x)x^{2}}{||\nu_{B^{c}}||+\int_\s \pi(\dd \theta)\mu^{d}(\theta)}\right|<\epsilon,$$ for all $t>T.$ We apply this in the following estimate, for $t>T,$ we have that there is a constant $K_1\in(0,\infty),$ such that
\begin{equation*}
\begin{split}
&\int_{T}^t  \frac{1}{(1+A_s)^2}\dd C_s + \int_T^t \dd s
\int_\s\Pi(\Theta_s, \s, \dd x)\frac{x^2}{1+A_s}\frac{1}{1+A_s + |x|}\one{|x|\leq 1}\\
&\leq \int_{T}^t  \frac{1}{(1+A_s)^2}\dd C_s +\int_T^t \dd s
\int_\s\Pi(\Theta_s, \s, \dd x)\frac{x^2}{(1+A_s)^{2}}\one{|x|\leq 1}\\
&= \int_{T}^t \frac{\dd (C_s+\int^s_{0}\dd u\int_\s\Pi(\Theta_u, \s, \dd x)x^2)}{\left(1+(C_s+\int^s_{0}\dd u\int_\s\Pi(\Theta_u, \s, \dd x)x^2\right)^{2}}\times\frac{\left(1+(C_s+\int^s_{0}\dd u\int_\s\Pi(\Theta_u, \s, \dd x)x^2\right)^{2}}{(1+A_s)^{2}}\\
&\leq K_1\times\frac{||\nu_{C}||+\int_{s}\pi(\dd\theta)\int_{|x|\leq 1}\Pi(\theta, \s, \dd x)x^{2}}{||\nu_{B^{c}}||+\int_\s \pi(\dd \theta)\mu^{d}(\theta)}\\
&\qquad \times \left(-\left(1+(C_s+\int^s_{0}\dd u\int_\s\Pi(\Theta_u, \s, \dd x)x^2\right)^{-1}\right)\Bigg|^{t}_{T}.
\end{split}
\end{equation*}
When $t\to\infty,$ the rightmost term converges a.s. to 
$$\left(\left(1+(C_T+\int^T_{0}\dd u\int_\s\Pi(\Theta_u, \s, \dd x)\right)^{-1}\right).$$
We estimate next the term related to the integral over the jumps of magnitude larger than $1.$ The argument we will use mimics the previous one, so we will skip some steps. By the assumption of having a $1+p$-finite moment, (H8), we obtain the upper bound
\begin{equation*}
\begin{split}
&\int_0^t \dd s
\int_\s\Pi(\Theta_s, \s, \dd x)\one{|x|>1}\frac{x^2}{1+A_s}\frac{1}{1+A_s + |x|}\\ 
&\leq \int_0^t \dd s
\int_\s\Pi(\Theta_s, \s, \dd x)\one{|x|>1}\frac{|x|^{1+p}}{(1+A_s)^{1+p}}\frac{|x|^{1-p}}{(1+A_s + |x|)^{1-p}}\\
&\leq \int_0^t \dd s
\int_\s\Pi(\Theta_s, \s, \dd x)\one{|x|>1}\frac{|x|^{1+p}}{(1+A_s)^{1+p}}.
\end{split}
\end{equation*}
Denoting by $R^{(p)}$, the additive functional 
\[
R^{(p)}_s=\int^s_0\dd u\Pi(\Theta_u, \s, \dd x)|x|^{1+p}\one{1<|x|<\infty},
\qquad s\geq0
\]
 we have by a further application of the Ergodic Theorem of X.3.12 of \cite{RY} that 
$$\frac{R^{(p)}_t}{A_t}\xrightarrow[t\to\infty]{\text a.s.}\frac{\int_{\s}\pi(\dd\theta)\int_{x\in\R}\Pi(\theta, \s, \dd x)|x|^{1+p}\one{1<|x|< \infty}}{m_1}.$$ As before, we have that there is $T$ and constant $K_2\in(0,\infty)$ such that for $t>T,$
\begin{equation*}
\begin{split}
&\int_T^t \dd s
\int_\s\Pi(\Theta_s, \s, \dd x)\one{1<|x|<\infty}|x|^{1+p}\left(\frac{1}{(1+A_s)^{1+p}}\right)\\
&\leq K_2 \int_T^t \dd R^{(p)}_s\left(\frac{1}{(1+R^{(p)}_s)^{1+p}}\right)\\
&\xrightarrow[t\to\infty]{\text{a.s.}} \frac{K_2}{p(1+R^{(p)}_T)^{p}}.
\end{split}
\end{equation*}

Putting the pieces together, we conclude that $\widetilde{G}_\infty<\infty,$ a.s. and hence that \eqref{tildaSLLNA} holds, from which the first claim in \eqref{tildaSLLN+}  follows. 
\medskip

For the second claim in  \eqref{tildaSLLN+}, note from  Theorem of X.3.12 of \cite{RY},  since $A$ is an additive functional, we have 
\begin{equation}
\frac{1}{t}A_t\xrightarrow[t\to\infty]{\text{a.s.}}\frac{m_1}{\pi(1)},
\label{At/t*}
\end{equation}
which equals zero as $\pi(1) = \infty$.

%When the process $\Theta$ is positive recurrent, by gathering \eqref{tildaSLLNA} and \eqref{At/t*}, we have the strong law of large numbers
%\begin{equation}\label{SLLN-Nfixed}
%\begin{split}
%\frac{\xi_t}{t}\xrightarrow[t\to\infty]{\text{a.s.}} ||\nu_{B^c}||&+\int_{\s}\pi(\dd \theta)\int_{\R\setminus\{0\}}x\one{1<|x|<\infty}\Pi(\theta, \s,\dd x)\\
%&+\int_{\s}\pi(\dd \theta)\int_{\R\setminus\{0\}}x\one{|x|\leq 1}\Pi(\theta, \s\setminus\{\theta\},\dd x).
%\end{split}
%\end{equation}
\medskip

Now suppose we are still in the setting that $\pi$ is infinite, but $m_1= 0.$ We can artificially add a linear drift with  rate $\kappa>0$ to $\xi$. Write $\xi^{(\kappa)}_t = \xi_t+ \kappa t$, $t\geq0$. We will similarly index other quantities such as $A^{(\kappa)}$ accordingly for the process $\xi^{(\kappa)}$. Note that $A^{(k)}_t = \kappa t + A_t$. The invariant measure for $\Theta$ remains the same.

\medskip

According to what we have proved above, we have almost surely
\begin{equation}
1= \lim_{t\to\infty}\frac{\xi^{(k)}_t}{A^{(k)}_t } =  \lim_{t\to\infty}\frac{\xi _t + \kappa t}{A_t + \kappa t} =  \lim_{t\to\infty}\frac{\xi _t} {t ( (A_t/t)+\kappa)} + 
\lim_{t\to\infty} \frac{1}{((A_t /t\kappa) + 1)},
\label{justifythelimits}
\end{equation}
where the first limit on the right-hand side is justified if the second  limit on the right-hand side is justified. However, this is the case, thanks to \eqref{At/t*}.
Back in \eqref{justifythelimits}, we can now see that 
\[
\lim_{t\to\infty}\frac{\xi_t}{t }  = 0
\]
almost surely.
\medskip

Finally we turn to the setting that $\pi(1) = 1$. Technically, we may repeat the arguments above, noting that in \eqref{At/t*} the limit is precisely $m_1 = \bE_{0,\pi}[\xi_1]$.
However, we can do better by revisiting Proposition 2.15 of \cite{TAMS}.
% use Lipster's theorem to prove directly that \eqref{tildaSLLNA} holds with $A_t$ replaced by $t.$ All the arguments provided after \eqref{tildaSLLNA} hold with small adaptations to ensure that in this case we also have $\widetilde{G}_{\infty}<\infty,$ a.s. 

%We denote now by $(\Theta, \xi^{N}),$ the MAP which is obtained by taking
%$$\xi^{N}_t=\sum_{s\leq t}\Delta \xi_s\one{|\Delta \xi_s|>N},\qquad t\geq 0.$$ $(\Theta, \xi^{N}),$ is a MAP of the pure jump type, as the ones considered by Breuer in \cite{breuer}. Actually, the jump kernel defined in (2) in page 77 in \cite{breuer} is given by the finite kernel $$\Pi(\theta, \dd \beta,\dd x)\one{|x|>N}, \qquad \theta,\beta\in\s, x\in\R.$$ 
%Our assumptions ((H2) and (H5), precisely) imply that the assumptions in Theorem 27 in \cite{breuer} are fulfilled. From there we deduce that  
%\begin{equation}\label{SLLN-Nfixedbis}
%\frac{\xi^N_t}{t}\xrightarrow[t\to\infty]{\text{a.s.}}\int_{\s}\pi(\dd \theta)\int_{\R\setminus\{0\}}x\one{N<|x|}\Pi(\theta, \s,\dd x).
%\end{equation}
%Adding the estimates obtained in \eqref{SLLN-Nfixed} and  \eqref{SLLN-Nfixedbis},  

Indeed, recalling that (H7) implies (H5), we can write 
% Assume that  \eqref{M1} is satisfied. The centered process 
\begin{equation}
\begin{split}
\widehat{\xi}_t&=\xi_t-
\bQ^{\Theta}[\xi_t]\\
&=\xi_t-B^{c}-\sum_{0<s\leq t}\int_{\R}y\nu_{(\Theta_{s-}, \Theta_{s})}(\dd y)\one{y\neq 0}\\
&\qquad -\int^{t}_{0}\dd u\int_{\R\setminus\{0\}}x\one{|x|>1}\Pi(\Theta_u, \{\Theta_u\},\dd x), \qquad t\geq 0,
\end{split}
\end{equation} is such that $( \widehat{\xi},\Theta)$ is a MAP. Moreover,  conditionally on $\Theta$ it has independent increments. 
Doob's Theorem 5.1 in page 363 in \cite{doobbook} can be applied to it to show that, under $\bQ^{\omega}$
\begin{equation*}
\bQ^{\omega}[\sup_{0\leq s\leq t}|\widehat{\xi}_s|]\leq 8 \bQ^{\omega}[|\widehat{\xi}_t|],\qquad t\geq0,
\end{equation*}
which implies that 
\begin{equation*}
\bE_{0,\theta}[\sup_{0\leq s\leq t}|\widehat{\xi}_s|]\leq 8 \bE_{0,\theta}[|\widehat{\xi}_t|],\qquad t\geq0.
\end{equation*}
If (H5) is satisfied, the triangle inequality implies that the rightmost term above is finite. This should be useful in proving the moment condition for the process $\xi.$ Indeed, we have by the triangle inequality that for any $t\geq 0$
\begin{equation}
\begin{split}
\bE_{0,\theta}[\sup_{0\leq s\leq t}|\xi_s|]&\leq \bE_{0,\theta}[\sup_{0\leq s\leq t}|\widehat{\xi}_s|]+\bE_{0,\theta}[\sup_{0\leq s\leq t}|\bQ^{\Theta}[\xi_t]|]\\ 
&\leq 8 \bE_{0,\theta}\left[|\widehat{\xi}_t| \right]+ \bE_{0,\theta}\left[\sup_{0\leq u\leq t}|B^{c}_u|\right]\\&+\bE_{0,\theta}\left[\int^{t}_{0}\dd u\int_{\R\setminus\{0\}}|x|\one{|x|>1}\Pi(\Theta_u, \s,\dd x)\right]\\
&+\bE_{0,\theta}\left[\int^{t}_{0}\dd u\int_{\R\setminus\{0\}}|x|\one{|x|\leq 1}\Pi(\Theta_u, \s\setminus\{\Theta_u\},\dd x)\right].
\end{split}
\end{equation}
Noting that $
\bE_{0,\theta}\left[\sup_{0\leq u\leq t}|B^{c}_u|\right]\leq  \bE_{0,\theta}\left[\int^{t}_{0}|b_s| \dd s\right]
$, the right-hand side above is finite thanks to the assumption (H7), which implies (H5).
\hfill$\square$

%\begin{proof}
%The proof follows as a direct consequence of the Theorem~3.16 in~\cite{HandL} once we notice the following Lemma.
%\end{proof}
%\begin{proof}[Proof of Lemma~(\ref{square})]According to (\eqref{xisemi-mart}), and Proposition 2.29 in Chapter II in \cite{J&S}, 
%\begin{equation*}
%(\xi_t-\widetilde{B}_t)^{2}-C_t-\int^{t}_{0}\dd s\int_{\R\setminus\{0\}} x^{2}\Pi(\Theta_s, \{\Theta_s\}, \dd x), \qquad t\geq 0,
%\end{equation*}
%is a local martingale under $\bQ^{\Theta}.$ Besides, arguing as in \eqref{display9}, we have that the dual projection of 
%$$\sum_{0<s\leq t}\int_{\R}y^{2}\nu_{(\Theta_{s-}, \Theta_{s})}(\dd y), \qquad t\geq 0,$$ is given by 
%$$\int^{t}_{0}\dd s\int_{\R}y^{2}\Pi(\Theta_{s}, \s\setminus\{\Theta_{s}\}, \dd y), \qquad t\geq 0.$$
%
%
%Hence, 
%\begin{equation*}
%\begin{split}
%\bE\left(\right)
%\end{split}
%\end{equation*}
%
%\end{proof}
\section{Proof of Theorem \ref{CLTMAPS}}
Our proof will predominantly use the main results from \cite{HandL}.
We assume that (H3), (H4) and (H6) hold, as well as (H8) with $p=1.$ The latter and Proposition 2.29 in~\cite{J&S} imply that $\xi$ is a square integrable semimartingale and also that the bracket process $$\langle\xi\rangle_t=C_t+\int^{t}_0\dd s\int_{\R\setminus\{0\}}\Pi(\Theta_s, \s, \dd y)y^{2},\qquad t\geq 0,$$ is an additive functional of $\Theta,$ and hence it is $\Theta$ measurable. To prove the Theorem we will follow the usual path determining the limit through finite dimensional convergence, and then we will prove tightness. 

\medskip
For finite dimensional convergence, since we have a process with conditionally independent increments, it will actually be enough to study the one dimensional distributions. This will be our first task.

As before, by the Ergodic Theorem of X.3.12 of \cite{RY} we have that
\begin{equation}
\frac{\langle\xi\rangle_t}{t}\xrightarrow[t\to\infty]{\text{a.s.}}\frac{||\nu_C||+\int_{\s}\pi(\dd\theta)\int_{x\in\R\setminus\{0\}}\Pi(\theta, \s, \dd x)x^{2}}{\pi(1)}.
\end{equation}
Recall the decomposition in \eqref{xisemi-mart}, and the notation in \eqref{Btilde}. We will consider  the characteristic function of $$\frac{\xi_{nt}-A_{nt}}{\sqrt{h(n)}}=\frac{\xi_{nt}-\widetilde{B}_{nt}}{\sqrt{h(n)}}+ \frac{\widetilde{B}_{nt}-A_{nt}}{\sqrt{h(n)}},  \qquad t\geq 0.$$
That is, 
\begin{equation*}
\begin{split}
\bE_{0,\theta}\left[\exp\left\{i\lambda \frac{\xi_{nt}-A_{nt}}{\sqrt{h(n)}}\right\}\right]=\bE_{0,\theta}\left[\exp\left\{i\lambda \frac{\widetilde{B}_{nt}-A_{nt}}{\sqrt{h(n)}}\right\}\bQ^{\Theta}\left[\exp\left\{i\lambda \frac{\xi_{nt}-\widetilde{B}_{nt}}{\sqrt{h(n)}}\right\}\right]\right].
\end{split}
\end{equation*}
Using equation \eqref{eq:charac} and some algebra, we obtain the expression
\begin{equation}
\begin{split}
&\bQ^{\Theta}\left[\exp\left\{i\lambda \frac{\xi_{nt}-\widetilde{B}_{nt}}{\sqrt{h(n)}}\right\}\right]\\
&=\exp\Bigg\{\Bigg(-\frac{\lambda^{2}}{2h(n)}C_{nt}\notag\\
&\hspace{3cm}+\int^{nt}_{0}\dd u\int_{\R\setminus\{0\}}\left(e^{i\lambda\frac{x}{\sqrt{h(n)}}}-1-i\lambda\frac{x}{\sqrt{h(n)}} \right)\Pi(\Theta_{u}, \{\Theta_{u}\},\dd x)\Bigg)\Bigg\}\\
& \times \prod_{0<s\leq nt}\exp\left\{-\int_{\R}\frac{i\lambda y}{\sqrt{h(n)}}\nu_{(\Theta_{s-}, \Theta_{s})}(\dd y)\one{y\neq 0}\right\}\notag\\
&\hspace{4cm}\cdot\left(1+\one{\Theta_{s-}\neq \Theta_s}\int_{\R}(e^{\frac{i x\lambda}{\sqrt{h(n)}}}-1)\nu_{(\Theta_{s-}, \Theta_s)}(\dd x)\right).
\end{split}
\end{equation}
Then, some elementary manipulations allow us to write  
\begin{equation}
\begin{split}
&\exp\left\{i\lambda \frac{\widetilde{B}_{nt}-A_{nt}}{\sqrt{h(n)}}\right\}\bQ^{\Theta}\left[\exp\left\{i\lambda \frac{\xi_{nt}-\widetilde{B}_{nt}}{\sqrt{h(n)}}\right\}\right]\\
&=M^{(1)}_t\times M^{(2)}\times M^{(3)}_t,\\
\end{split}
\end{equation}
where 
\begin{equation}
\begin{split}
M^{(1)}_t&:=\exp\Bigg\{\Bigg(-\frac{\lambda^{2}}{2h(n)}C_{nt}\notag\\
&\hspace{2cm}+\int^{nt}_{0}\dd u\int_{\R\setminus\{0\}}\left(e^{i\lambda\frac{x}{\sqrt{h(n)}}}-1-i\lambda\frac{x}{\sqrt{h(n)}} \right)\Pi(\Theta_{u}, \{\Theta_{u}\},\dd x)\Bigg)\Bigg\}\\
&\times\exp\left\{\sum_{s\leq nt}\int_{\R\setminus\{0\}}\nu_{(\Theta_{s-}, \Theta_{s})}(dx)\left(e^{i\lambda\frac{x}{\sqrt{h(n)}}}-1-i\lambda\frac{x}{\sqrt{h(n)}} \right)\right\}, \qquad t\geq 0;
\end{split}
\end{equation}
\begin{equation}
\begin{split}
M^{(2)}_{t}:=&\prod_{0<s\leq nt}\left(1+\int_{\R}\left(e^{\frac{i x\lambda}{\sqrt{h(n)}}}-1\right)\nu_{(\Theta_{s-}, \Theta_s)}(\dd x)\right)\\
&\times\prod_{0<s\leq nt}\exp\left(-\int_{\R}\left(e^{\frac{i x\lambda}{\sqrt{h(n)}}}-1\right)\nu_{(\Theta_{s-}, \Theta_s)}(\dd x)\right), \qquad t\geq 0;
\end{split}
\end{equation}
and
\begin{equation}
\begin{split}
M^{(3)}_{t}&:=\exp\left\{\frac{i\lambda}{\sqrt{h(n)}}Z_{nt}\right\},\qquad t\geq 0,
\end{split}
\end{equation}
where 
\begin{equation}\label{Zt}
Z_{nt}:=\sum_{s\leq {nt}}\int_{\R\setminus\{0\}}\nu_{(\Theta_{s-}, \Theta_{s})}(dx)x-\int^{nt}_{0}\dd u\int_{\R\setminus\{0\}}\Pi(\Theta_{u}, \s\setminus\{\Theta_u\},\dd x)x,\qquad t\geq 0.
\end{equation}

For  $M^{(1)}$ and $M^{(2)}$, we recall that there are $a_1, a_2\in\mathbb{C},$ such that $|a_1(\cdot)|<1,$ $|a_2(\cdot)|<1$ and
\begin{equation*}
 e^{iy}-1-iy= a_1(y)\frac{y^{2}}{2}, \qquad e^{iy}-1-iy-\frac{(iy)^{2}}{2}= a_2(y) \frac{|y|^{3}}{6},\qquad \text{for any}, \ y\in\R;
\end{equation*}
and $$\left|e^{iy}-1-iy\right|=O(y^{2}), \qquad y>1;$$
see e.g. Lemma 8.6 in \cite{Satobook}. Recall the random measure $L$ defined in \eqref{L}. Using this notation, and the above facts, we write $M^{(1)}$ as
\begin{equation}
\begin{split}
M^{(1)}_t&=\exp\left\{-\frac{\lambda^{2}C_{nt}}{h(n)}+\int_{(0,nt]\times\R\setminus\{0\}}L(\dd s, \dd x)\left(e^{\frac{i\lambda x}{\sqrt{h(n)}}}-1-\frac{i\lambda x}{\sqrt{h(n)}}\right)\right\}.
\end{split}
\end{equation}

In order to deal with $M^{(2)}$ and $M^{(3)}$, we must first make some additional observations.  First, by the  Ergodic Theorem of X.3.12 of \cite{RY}, and a calculation similar to the  one in \eqref{display9}, we obtain that for any $\epsilon>0$
\begin{equation}\label{variance1}
\begin{split}
&\frac{1}{\langle\xi\rangle_{nt}}\left(-\lambda^{2}C_{nt}+\int_{(0,nt]\times\R\setminus\{0\}}L(\dd s, \dd x)\left(e^{\frac{i\lambda x}{\sqrt{h(n)}}}-1-\frac{i\lambda x}{\sqrt{h(n)}}\right)\one{|x\lambda|\leq\epsilon\sqrt{h(n)}}\right)\\
&\stackrel{t\to\infty}{\sim} \left(\frac{-\lambda^{2}C_{nt}-\lambda^{2}\int_{(0,nt]}L(\dd s, \dd x)x^{2}\one{|\lambda x|\leq \epsilon\sqrt{h(n)}}}{\langle\xi\rangle_{nt}}\right)\xrightarrow[t\to\infty]{\text{a.s.}} -\frac{\lambda^{2}}{2}.
%\left(\frac{||\nu_C||+\int_{\theta\in\Theta}\pi(\dd \theta)\int_{x\in\R}\Pi(\theta, \s, \dd x)x^{2}}{||\nu_C||+\int_{\theta\in\Theta}\pi(\dd \theta)\int_{x\in\R}\Pi(\theta, \s, \dd x)x^{2}}\right).
\end{split}
\end{equation}

Second, due to the Lindenberg condition \eqref{lin}, we also  have
\begin{equation}\label{variance1.2}
\begin{split}
&\frac{1}{h(nt)}\left(\int_{(0,nt]\times\R\setminus\{0\}}L(\dd s, \dd x)\left(e^{\frac{i\lambda x}{\sqrt{h(n)}}}-1-\frac{i\lambda x}{\sqrt{h(n)}}\right)\one{|x\lambda|>\epsilon\sqrt{h(n)}}\right)\\
&\stackrel{}{\sim} \left(\frac{-\lambda^{2}\int_{(0,nt]}L(\dd s, \dd x)x^{2}\one{|\lambda x|> \epsilon\sqrt{h(n)}}}{h(n)}\right) \xrightarrow[t\to\infty]{\text{a.s.}} 0.
\end{split}
\end{equation}

Third, the process $Z$ defined in \eqref{Zt} is a square integrable martingale, whose angle bracket $\langle Z\rangle$ and square bracket $[Z]$ processes are additive functionals of $\Theta.$ (Recall that $[Z]$ is the  compensator of $Z^2$  and $\langle Z\rangle$ is the predictable compensator of $[Z]$.) Indeed, we can use the compensation formula for L\'evy systems, Jensen's inequality, and (H8), to write 
$$[Z]_t=\sum_{s\leq t}\left(\int_{\R\setminus\{0\}}\nu_{(\Theta_{s-}, \Theta_{s})}(dx)x\right)^{2},\qquad t\geq 0;$$
$$\langle Z\rangle_t=\int^{t}_{0}\dd s \int_{\s}K(\Theta_s, \dd \beta)\left(\int_{\R}\nu_{(\Theta_s, \beta)}(\dd x)x\right)^{2},\qquad t\geq 0;$$ and 

\begin{equation*}
\begin{split}
\bE_{0,\theta}\left[[Z]_t\right]=&\bE_{0,\theta}\left[\int^{t}_{0}\dd s \int_{\s}K(\Theta_s, \dd \beta)\left(\int_{\R}\nu_{(\Theta_s, \beta)}(\dd x)x\right)^{2}\right]\\
&\leq \bE_{0,\theta}\left[\int^{t}_{0}\dd s \int_{\s}K(\Theta_s, \dd \beta)\int_{\R}\nu_{(\Theta_s, \beta)}(\dd x)x^{2}\right]\\
&= \bE_{0,\theta}\left[\int^{t}_{0}\dd s \int_{\R}\Pi(\Theta_s, \s\setminus\{\Theta_s\}, \dd x)x^{2}\right]<\infty.
\end{split}
\end{equation*}
The Theorem 3.16 in \cite{HandL} implies that we have the joint convergence in the sense of Skorohod's topology
\begin{equation}\label{ML-convergence}
\begin{split}
\left(\left(\frac{1}{\sqrt{h(n)}}Z_{tn}, \,\frac{1}{h(n)}\langle Z\rangle_{nt}\right),\ t\geq 0\right)\xrightarrow[n\to\infty]{\text{Law}}\left(\left(J^{1/2}B_{W^{(\alpha)}_{t}}, \,J W^{(\alpha)}_{t}\right),\ t\geq 0\right);
\end{split}
\end{equation}
where $J$ is as defined in \eqref{J}. Also, by the Ergodic Theorem of X.3.12 of \cite{RY} and Slutsky's Theorem we have that we can extend the above joint convergence to  
\begin{equation}\label{triplet}
\begin{split}
&\left(\left(\frac{1}{\sqrt{h(n)}}Z_{tn}, \,\frac{1}{h(n)}\langle Z\rangle_{nt}, \,\frac{\langle\xi\rangle_{nt}}{h(n)}\right),\ t\geq 0\right)\\
&=\left(\left(\frac{1}{\sqrt{h(n)}}Z_{tn}, \,\frac{1}{h(n)}\langle Z\rangle_{nt},\, \frac{\langle Z\rangle_{nt}}{h(n)}\frac{\langle\xi\rangle_{nt}}{\langle Z\rangle_{nt}}, \right),\ t\geq 0\right)\\
&\xrightarrow[n\to\infty]{\text{Law}}\left(\left(J^{1/2}B_{W^{(\alpha)}_{t}}, \,J W^{(\alpha)}_{t}, \, J\times\frac{||\nu_{\langle\xi\rangle}||}{J} W^{(\alpha)}_{t}\right),\ t\geq 0\right).
\end{split}
\end{equation}

Now, to determine the limit of $M^{(1)}$,
we use \eqref{variance1} and the third component of the limit in \eqref{triplet}, together with \eqref{variance1.2}.
\medskip

To determine the limit of $M^{(2)}$ we recall that Taylor's representation of the logarithm allows to write
$$\log(1-z)+z=-\frac{z^{2}}{2}+O(|z|^{3}), \qquad |z|<1.$$ 
Hence, arguing as in \eqref{variance1} and \eqref{variance1.2}, appealing to the second element of \eqref{triplet}, we can see that 
\begin{equation*}
\begin{split}
M^{(2)}_{t}&=\exp\Bigg\{\sum_{s\leq nt}\Bigg[\log\left(1-\int_{\R}\left(1-e^{\frac{i\lambda x}{\sqrt{h(n)}}}\right)\nu_{(\Theta_{s-}, \Theta_s)}(\dd x)\right)\notag\\
&\hspace{5cm}+\int_{\R}\left(1-e^{\frac{i\lambda x}{\sqrt{h(n)}}}\right)\nu_{(\Theta_{s-}, \Theta_s)}(\dd x)\Bigg]\Bigg\}\\
&\sim \exp\left\{-\frac{\langle Z\rangle_{nt}}{2 h(n)}\frac{(-\lambda^{2})}{\langle Z\rangle_{nt}}\sum_{s\leq nt}\left[\int_{\R}\left(\frac{i\lambda x}{\sqrt{h(n)}}\right)^{-1}\left(1-e^{\frac{i\lambda x}{\sqrt{h(n)}}}\right)x\nu_{(\Theta_{s-}, \Theta_s)}(\dd x)\right]^{2}\right\}\\
&\xrightarrow[n\to\infty]{\text{a.s.}}\exp\left\{-\frac{\lambda^{2}J}{2}W^{(\alpha)}_{t}\right\}.
\end{split}
\end{equation*}
 Finally, to determine the limit of $M^{(3)}$, we can use the first element of \eqref{triplet}.
\medskip

Assembling the pieces we derive the limit
\begin{equation*}
\begin{split}
&\lim_{n\to\infty}\bE_{0,\theta}\left[\exp\left\{i\lambda \frac{\xi_{nt}-A_{nt}}{\sqrt{h(n)}}\right\}\right]\notag\\
&=\bE_{0,\pi}\left[\exp\left\{-\frac{\lambda^{2}}{2}\left(||\nu_{\langle\xi\rangle}||+J\right)W^{(\alpha)}_{t}+i\lambda J^{1/2}B_{W^{(\alpha)}_{t}}\right\}\right]\\
&=\bE_{0,\pi}\left[\exp\left\{i\lambda\left(||\nu_{\langle\xi\rangle}||+J\right)^{1/2}\Sigma^{(1)}_{W^{(\alpha)}_{t}}+i\lambda J^{1/2}\Sigma^{(2)}_{W^{(\alpha)}_{t}}\right\}\right];
\end{split}
\end{equation*}
where $\Sigma^{(1)}$ and $\Sigma^{(2)}$ are two independent standard Brownian motions, which are  independent of the Mittag-Leffler process $W^{(\alpha)}$.  The identification of the limit follows by conditioning on $W,$ and adding the variances of the two Gaussian processes.
\medskip

Tightness follows from the Theorem 4.13 in Chapter VI Section 4 in\cite{J&S}. Indeed, according to that result, it is enough to prove that the quadratic variation process of the martingale 
\[
\frac{1}{\sqrt{h(n)}}(\xi_{nt}-A_{nt}), \qquad t\geq 0,
\]
 is $C$-tight. That is,  every convergent subsequence converges to an a.s. continuous process. But this is a direct consequence of the convergence in \eqref{ML-convergence}, since the Mittag-Lefler process $W^{(\alpha)}$ is non-decreasing and continuous almost surely.
\hfill$\square$

\bibliography{references}{}
\bibliographystyle{plain}

\end{document}